\documentclass{amsart}

\usepackage{amsthm}

\usepackage{amsfonts}
\usepackage{graphicx}
\usepackage{amsmath}
\usepackage{amssymb}
\DeclareGraphicsExtensions{.pdf,.png,.jpg} %solo para PDFLaTeX
\usepackage{epstopdf}

\newtheorem{teorema}{Theorem}[section]

\newtheorem{proposicion}[teorema]{Proposition}
\newtheorem{pregunta}[teorema]{Question}

\theoremstyle{definition}

\theoremstyle{definition}
\newtheorem{definicion}[teorema]{Definition}

\theoremstyle{definition}

\newcommand{\ra}{\rightarrow}
\newcommand{\pr}{\partial}
\newcommand{\txt}{\textnormal}
\newcommand{\tl}{\tilde}
\newcommand{\empt}{\emptyset}

\newcommand{\minus}{\setminus}

\newcommand{\R}{\mathbf R}

\newcommand{\N}{\mathbf N}
\newcommand{\T}{\mathbf T}
\newcommand{\Z}{\mathbf Z}

\newcommand{\cl}{\mathcal}

\title{Dynamical models for some torus homeomorphisms}

\author{P. D\'avalos}

\date{}

\begin{document}

\maketitle

\begin{abstract}
Consider a homeomorphism of the torus $\mathbf{T}^2$ in the homotopy class of the identity. There is a topological invariant for $f$ known as the \textit{rotation set}, which is a compact convex subset of $H_1(\mathbf{T}^2,\mathbf{R}) \simeq \mathbf{R}^2$ and describes the homological direction and speed at which the orbits rotate on $\mathbf{T}^2$. In this paper we collect some results about the existence of dynamical models associated to this invariant. 
\end{abstract}

\section{Introduction}           \label{sec.intro}

In \cite{po} H. Poincar\'e defined the \textit{rotation number} for circle homeomorphisms and proved it to be a topological invariant carrying dynamical information. For an orientation preserving homeomorphism $f:\T^1\ra\T^1$ and a lift $\tl{f}:\R \ra\R$ of $f$, he defined the rotation number of $\tl{f}$ as the limit $\rho(\tl{f}):=\lim_{n\ra\infty} (\tl{f}^n(x)-x)/n$, which he showed to exist and to be independent of $x\in\R$. From the topological point of view, the dynamics of $f$ can be completely understood from the rotation number $\rho(\tl{f})$. This number is rational if an only if there exist periodic orbits for $f$, all of the same period, and the complement of the periodic points of $f$ is a union of periodic open intervals also with the same period. The rotation number $\rho(\tl{f})$ is irrational if and only if $f$ is semiconjugated to the rotation of $\T^1$ by the angle $\rho(\tl{f})$. If the semiconjugacy is not an actual conjugacy, then it is injective on a cantor set $K\subset\T^1$, and colapses the connected components of $\T^1\minus K$ to points. The set $K$ is minimal for $f$, and the connected components of $\T^1\minus K$ are wandering open intervals.

Later, in \cite{mz} Misiurewicz and Ziemian generalized the concept of the rotation number for homeomorphisms of $\T^n$. For a homeomorphism $f:\T^n\ra\T^n$ in the homotopy class of the identity, and a lift $\tl{f}:\R^n\ra\R^n$, the \textit{rotation set} of $\tl{f}$, denoted $\rho(\tl{f})$, is defined as the set of accumulation points of sequences of the form
$$ \left\{ \frac{\tl{f}^{n_i}(x)-x}{n_i} \right\}_{i\in\N}$$
with $x\in\R^2$ and $n_i\ra\infty$ as $i\ra\infty$. We will work with case that $n=2$. In this case the rotation set $\rho(\tl{f})$ is a compact convex subset of $\R^2$ \cite{mz} and it is also a topological invariant carrying dynamical information. For example, rational points of $\rho(\tl{f})$ are related to periodic orbits of $f$ \cite{f2, f3} and ergodic measures of $f$ are realated to extremal points of $\rho(\tl{f})$ \cite{mz}. 

Trying to emulate the case of $\T^1$, one could ask if there exist dynamical models associated to the rotation set. That is; to what extent, from the topological point of view, can one classify the dynamics of a torus homeomorphism from its rotation set? In this paper we gather some results related to this question, in the case that the rotation set is either a point or a segment.

\section{Notations and preliminaries}

Consider a curve $\gamma:I\ra\T^2$. We will denote its image Im$(\gamma)\subset\T^2$ also by $\gamma$, and we will call it also a \textit{curve}. We say that a curve $\gamma\subset\T^2$ is \textit{essential} if it is not homotopically trivial, and we say that $\gamma$ is \textit{vertical} if it is homotopic to a vertical circle $\{x\}\times\T^1\subset\T^2$. Similarily, we say that a topological annulus $A\subset\T^2$ is \textit{essential} if it is not homotopic to a point, and we say that $A$ is \textit{essential and vertical} if it is homotopic to a vertical annulus $ I \times \T^1$, with $I\subset\T^1$ an interval. A set $K\subset\T^2$ is said to be \textit{annular} if $K=\cap_{n\geq 0} A_i$, where the $A_i$ are topological compact annuli with $A_{i+1}\subset A_i$ such that $A_{i+1}\hookrightarrow A_i$ is a homotopy equivalence. If the annuli $A_i$ are essential, the annular set $K$ is called \textit{essential}, and if the $A_i$ are essential and vertical, the set $K$ is called \textit{essential and vertical}. We say that a set $K\subset\T^2$ is \textit{fully essential} if its complement is a union of pairwise disjoint open topological discs.

A homeomorphism $f:\T^2\ra\T^2$ is called a \textit{pseudo-rotation} if for a lift $\tl{f}$ (and hence for any lift), the limit
\begin{equation}    \label{eq1}
\lim_{n\ra\infty} \frac{\tl{f}^n(x)-x}{n}
\end{equation}
exists and is independent of $x$. The limit (\ref{eq1}) is called the \textit{rotation vector} of $\tl{f}$, and in the case that the rotation vector is zero, we say that $f$ and its lift $\tl{f}$ are \textit{irrotational}. 

A pseudo-rotation $f$ is said to have the \textit{bounded mean motion property} if there is a lift $\tl{f}$ with rotation vector $v$, such that the deviations
$$ D(x,n) = | \tl{f}^n(x)-x - nv | $$
are uniformly bounded in $x$ and $n$.

In general, the limit (\ref{eq1}) does not necessarily exist for every $x\in\R^2$, and we have the following definition.

\begin{definicion}[\cite{mz}]     %\label{def.rotset}
The \textit{rotation set of $\tl{f}$} is defined as 
$$\rho(\tl{f})=\bigcap_{m=1}^{\infty} \txt{cl} \left( \bigcup_{n=m}^{\infty} \left\{ \frac{\tl{f}^n(x)-x}{n}\, : \, x\in\R^2 \right\} \right) \subset\R^2. $$
%The \textit{rotation set of a point} $x\in\R^2$ is defined by
%$$\rho(x,f)= \bigcap _{m=1}^{\infty} \txt{cl} \left\{ \frac{f^n(x)-x}{n}\, : \, n>m \right\}.$$
%If the above set consists of a single point $v\in\R^2$, we call $v$ the \textit{rotation vector of} $x$. 
If $\Lambda\subset\T^2$ is a compact $f$-invariant set, we define the \textit{rotation set of} $\Lambda$ as
$$\rho(\Lambda,\tl{f})=\bigcap_{m=1}^{\infty} \txt{cl} \left( \bigcup_{n=m}^{\infty} \left\{ \frac{\tl{f}^n(x)-x}{n}\, : \, x\in\pi^{-1}(\Lambda) \right\} \right) \subset\R^2,$$
where $\pi:\R^2\ra\T^2$ denotes the canonical projection.
\end{definicion}

%\begin{remark}       %     \label{r.pot}
It is easy to see that for integers $n,m_1,m_2$,
$$\rho(T_1^{m_1}T_2^{m_2} \tl{f}^n)= n\rho(\tl{f}) + (m_1,m_2).$$
Then, the rotation set of any other lift of $f$ is an integer translate of $\rho(\tl{f})$, and we can think of the `rotation set of $f$' defined modulo $\Z^2$. 
%\end{remark}

\begin{teorema}[\cite{mz}]    % \label{mz.cc}
Let $f:\T^2\ra\T^2$ be a homeomorphism, and let $\tl{f}:\R^2\ra\R^2$ be a lift of $f$. Then the rotation set set $\rho(\tl{f})$ is compact and convex.
\end{teorema}

Consider a homeomorphism $f:\T^2\ra\T^2$. We say that $f$ has \textit{disc-type dynamics} if there exists a topological open disc $U\subset\T^2$ such that $\T^2\minus U \subset \txt{Fix}(f)$. Now, $f$ is said to have \textit{weakly annular dynamics} if for any lift $\tl{f}$ of $f$ there is $M>0$ and $v\in\Z^2\minus\{0\}$ such that 
$$|\langle   \tl{f}^n(x)-x , v  \rangle| \leq M$$	 
for all $x\in\R^2$ and $n\in \Z$. Finally, we say that $f$ has \textit{annular dynamics} if there exists an annular essential set set $A$ which is periodic for $f$. In this case, if $q$ is the period of $A$, to understand the dynamics of $f$ one may study the dynamics of the restrictions of $f^q$ to the open annulus $\T^2\minus A$, and to the annular set $A$. 

If the dynamics of $f$ is annular, it is easy to see that the rotation set of any lift is a (possibly degenerate) interval with rational slope containing rational points. If the rotation set of $\tl{f}$ is a non-degenerate interval with rational slope containing rational points, we will say that $\tl{f}$ has \textit{annular rotation set} (cf. Question \ref{preg1}).

Consider an orientable surface $M$ without boundary, and an isotopy $I:M \times [0,1]\ra M$ from the identity to a homeomorphism $f$. For a point $x\in M$, denote $\gamma_x= I(x,\cdot)$. A fixed point $p$ of $f$ is called \textit{contractible} if the closed path $\gamma_p$ is not homotopically trivial. An oriented topological foliation $\cl{F}$ of $M$ is said to be \textit{transverse} to $I$ if for any $x\in M$, the path $\gamma_x$ is homotopic with fixed endpoints to a path which is positively transverse to $\cl{F}$. Having contractible fixed points is clearly an obstruction to the existence of a foliation of $M$ transverse to $I$. The following theorem from \cite{lc2} says that it is actually the only obstruction. 

\begin{teorema}   \label{teo.pat}
If $f$ has no contractible fixed points, then there exists a topological oriented foliation without singularities which is transverse to the isotopy $I$. 
\end{teorema}

In a similar way, if $\cl{F}$ is a foliation with singularities, we say that $\cl{F}$ is transverse to the isotopy $I$ if for any $x\in M\minus\txt{sing}(\cl{F})$, the path $\gamma_x$ is homotopic with fixed endpoints to a path which is positively transverse to $\cl{F}\minus\txt{sing}(\cl{F})$. The following result is a consequence of \cite{lc2} and \cite{jau}. 

\begin{teorema}    \label{jaulent}
If $\txt{Fix}(f)$ is totally disconnected, then there exists a compact set $X\subset\txt{Fix}(f)$, an oriented foliation $\cl{F}$ with singularities in $X$, and an isotopy $I$ from the identity to $f$ such that $I$ fixes $X$, and $\cl{F}$ is transverse to $I$.
\end{teorema}

\section{Pseudorotations}

The canonical model of a pseudo-rotation is a rigid rotation. There exist both \textit{positive} and \textit{negative} results, in the sense that certain hypotheses guarantee or not some similarity of the dynamics of a pseudo-rotation and the corresponding rigid rotation.

One way to relate a pseudo-rotation and a rotation is by means of a semi-conjugacy. The question of whether such a semiconjugacy exists or not, has been systematically studied by T. J\"ager. We start by stating a basic result of this kind. From now on, for $\rho\in\R^2$, $R_{\rho}$ will denote the rotation $x\mapsto x+\rho \, \mod \Z^2$ on $\T^2$.

\begin{proposicion}      \label{jager1}
Let $f\in\txt{Homeo}_0(\T^2)$ be a minimal pseudo-rotation with bounded mean motion and a totally irrational rotation vector $\rho\in\R^2$. Then $f$ is semiconjugate to $R_{\rho}$. 
\end{proposicion}

In \cite{j1} it is proven a more general version of this proposition, for pseudorotations of $\T^n$, $n\geq 2$ (result which also deals with rotation sets that are not reduced to a point). 

The proof of Proposition \ref{jager1} consists in obtaining a semiconjugacy $h_i:\T^2\ra\T^1$ of the map $f$ with the rotation in $\T^1$ by $\rho_i$, for $i=1,2$. Then, the semiconjugacy $h:\T^2\ra\T^2$ between $f$ and $R_{\rho}$ is defined as $h=(h_1,h_2)$. The $h_i$ are defined as 
$$h_i(z) = \sup_{n\in\Z} (\pi_i \circ F^n(z) - n \rho_i),$$
for $i=1,2$. Due to the bounded mean motion property, the $h_i$ are well defined, and it is easy to check that $H(z)+ \rho_i = H \circ F(z)$. The minimality of $f$ is then used to prove continuity, and surjectivity is due to the minimality of the one dimensional rotation by $\rho_i$. 

For conservative pseudo-rotations, we have the following classification result.

\begin{teorema}[\cite{j1}]    \label{jager2}
Suppose $f \in \txt{Homeo}_0(\T^2)$ is a conservative pseudo-rotation with rotation vector $\rho\in\R^2$ and bounded mean motion. Then the following hold:
\begin{enumerate}
\item $\rho$ is totally irrational if and only if $f$ is semi-conjugate to $R_{\rho}$.
\item $\rho$ is neither totally irrational nor rational if and only if $f$ has a periodic circloid.
\item $\rho$ is rational if and only if $f$ has a periodic point.
\end{enumerate}
\end{teorema}

Item 3 of this theorem is a classical result of Franks \cite{f1}. For item (2), suppose that $\rho$ is a horizontal vector $(a,0)$, with $a$ irrational. In this case J\"ager proves that the bounded mean motion property is equivalent to having an actual horizontal `barrier' for the dynamics. Precisely, he proves that for such $\rho$, there exists a horizontal $f$-invariant circloid if and only if $f$ has the bounded mean motion property (a horizontal circloid is a circloid contained in an annulus homotopic to the horizontal annulus $\T^1\times [0,1/2]$). Then, the case for a general $\rho$ that is neither totally irrational nor rational is then deduced from this. 

%in \cite{j1} J\"ager studies periodic circloids for non-wandering torus homeomorphisms, and he proves that if $\rho$ is a vertical vector $(a,0)$, with $a$ irrational, then the bounded mean motion property is equivalent to having an actual horizontal `barrier' for the dynamics; that is, $f$ has the bounded mean motion property if and only if there exists an $f$-invariant horizontal circloid (a horizontal circloid is a circloid contained in an annulus homotopic to the horizontal annulus $\T^1\times [0,1/2]$). The case for a general $\rho$ that is neither totally irrational nor rational is then deduced from this. 

In item (1), the `if' part is elementary. The `only if' part of this item is of principal interset in the paper. One works with the lift $\hat{f}$ of $f$ to $\T^1\times\R$ with average vertical displacement $\rho_2$, and such that $|\pi_2 (\hat{f}^n(z) - z) - n\rho_2 | \leq c \ \forall n\in\Z, \, z\in\T^1\times \R$, where $c$ is the bounded mean motion constant for $f$. 

For $r\in\R$, one defines the sets
$$A_r = \bigcup_{n\in\Z} \hat{f}^n(\T^1 \times \{ r - n\rho_2 \} ).$$
By the mean motion property, the sets $A_r$ are bounded in the vertical direction. A main step here carried out in \cite{j1} is to extract from each $A_r$ a circloid $C_r$, in a way that the circloids are pairwise disjoint, and if $T :\T^1\times \R \ra \T^1\times\R$ denotes the translation $T(x,y)=(x,y+1)$, then
\begin{equation}  \label{eq.circ1}
C_{r+1} = T(C_r),
\end{equation}
\begin{equation}   \label{eq.circ2}
\hat{f}(C_r) = C_{r+\rho_2},
\end{equation}
and
\begin{equation}    \label{eq.circ3}
C_r\prec C_s \txt{ if } r < s,
\end{equation}
where the notation $C_r \prec A$ means that $A$ is contained in the connected component of $(\T^1\times \R) \setminus C_r$ that is unbounded from above and bounded from below. 

Having showed this, one can construct a semiconjugacy $H_2:\T^2\ra\T^1$ between $f$ and the one dimensional rotation by $\rho_2$ in the following way. Set
$$H_2(z) = \sup \{ r\in\R \, | \, C_r \prec z \}.$$
Using (\ref{eq.circ1}) and (\ref{eq.circ2}) it can be easily verified that
\begin{equation}   \label{eq.semic1}
H_2\circ T(z) = H(z) + 1,
\end{equation} 
and
\begin{equation}    \label{eq.semic2}
H_2\circ \hat{f}(z) = H(z) + \rho_2.
\end{equation}
To see that $H_2$ is continuous, it is proved that for any open interval $(a,b)\subset\R$, the set $H_2^{-1}(a,b)$ is a union of open `intervals' of the form $(C_r,C_s)\subset \T^1\times \R$, $r<s$. Here, $(C_r,C_s) = \{ z\in\T^1\times\R \, | \, C_r \prec z \prec C_s\}$. By properties (\ref{eq.semic1}) and (\ref{eq.semic2}), $H_2$ projects then to a semiconjucacy $h_2$ between $f$ and the irrational rotation $x\mapsto x + \rho_2$. In the same way, one constructs a semiconjugacy $h_1$ between $f$ and the rotation $x\mapsto x + \rho_1$, and then $h=(h_1,h_2)$ gives a semiconjugacy between $f$ and $R_\rho$ on $\T^2$.

The conservative hypothesis in Theorem \ref{jager2} is necessary, as the following proposition shows.

\begin{proposicion}[\cite{j1}]
Given any totally irrational rotation vector $\rho\in\R^2$, there exists an irrational pseudo-rotation $f\in\txt{Homeo}_0(\T^2)$ which has rotation vector $\rho$ and bounded mean motion, but which is not semi-conjugate to the irrational rotation $R_{\rho}$.
\end{proposicion}

One can easily see that the bounded mean motion property is also necessary in order to have a semiconjugacy between $f$ a rotation $R_\rho$. Actually, the absence of this property allows to create exotic examples of pseudo-rotations, with dynamical properties far from the rigid rotations. An illustration of this phenomenon is given by a result of Koropecki and Tal, for the irrotational case. Before stating that result, we give a definition. Given a set $X\subset\R^2$, we say that $X$ accumulates at infinity in the direction $v\in\T^1$ if there is a sequence of $x_n\in X$ such that $|x_n|\ra\infty$ and $(x_n-x_0)/|x_n-x_0|\ra v$ as $n\ra\infty$. The boundary of $X$ at infinity is defined as the set $\pr_{\infty}X$ consisting of all $v\in\T^1$ such that $X$ accumulates in the direction $v$ at infinity.

\begin{teorema}[\cite{korotal}]            \label{teo.kt1}
There exists a $C^{\infty}$ area-preserving irrotational pseudo-rotation which is ergodic with respect to Lebesgue measure, and such that for almost every point $x\in\R^2$, the orbit of $x$ accumulates in every direction at infinity, i. e., 
$$\pr_{\infty} \{\tl{f}^n(x)\, : \, n\in\N\} = \T^1.$$
\end{teorema}

For the case that the rotation vector $\rho$ is any rational vector, an analogous statement is obtained for a power of $f$. The proof of Theorem \ref{teo.kt1} theorem uses the idea of embedding an open disc $U$ in $\T^2$ in a way that $U$ has full Lebesgue-measure and any lift of $U$ to the universal covering accumulates in every direction at infinity. Then, one `glues' in $U$ a diffeomorphism of the open disc which is ergodic with respect to Lebesgue (or even isomorphic to a Bernoulli shift) and extends to the boundary of $U$ as the identity. The construction of such a diffeomorphism of the unit disc can be done by slight modifications of results due to Katok \cite{ka1}.

More examples of pseudo-rotations without the bounded mean motion property and with dynamical properties distant from the rigid rotation, both from the topological and metrical point of view, can be found for example in \cite{j2}, \cite{kk2} and \cite{fay1}.

Koropecki and Tal have also found that the `exotic' construction from Theorem \ref{teo.kt1} is to some extent the only way to create this behavior. That is, the existence of unbounded mean motion in many directions forces the existence of a `large' (fully essential) set of fixed points. This is shown in the following result, which classifies conservative irrational pseudo-rotations.

\begin{teorema}[\cite{korotal2}]
Let $f\in\txt{Homeo}_0(\T^2)$ be a conservative irrotational pseudo-rotation, and let $\tl{f}$ be its irrotational lift. Then one of the following holds (cf. Section \ref{sec.intro} for definitions):
\begin{enumerate}
\item $f$ has disc-type dynamics, and $\txt{Fix}(f)$ is fully essential,
\item every point in $\R^2$ has a bounded $\tl{f}$-orbit,
\item $f$ has weakly annular dynamics. 
\end{enumerate}
\end{teorema}

The proof of this theorem is done by contradiction. If the theorem does not hold, then Fix$(f)$ is not fully essential, $f$ is not annular, and there exists $x\in\R^2$ with unbounded $\tl{f}$-orbit. If $f$ is not annular, then one may show that Fix$(f)$ is inessential, and in \cite{korotal2} it is showed that one may actually assume that Fix$(f)$ is totally disconnected. A remarkable theorem proved in the same article (Theorem F) shows that the fact that $f$ is conservative and irrotational implies that the irrotational lift $\tl{f}$ is non-wandering. Now, the fact that $\tl{f}$ is non-wandering allows Koropecki and Tal to prove that the Brouwer foliation of $\R^2$ given by Proposition \ref{jaulent} is actually a \textit{gradient-like} foliation; that is, every leaf $\gamma$ is such that $\alpha(\gamma)=\{p_1\}$ and $\omega(\gamma)=\{p_2\}$, for some $p_1,p_2\in\tl{X}$, $p_1\neq p_2$. From there, using techniques from \cite{azt} and developing more machinery, it is worked to find a contradiction, which will prove the theorem.

\section{Rotation sets which are intervals}

The only known examples of rotation sets which are intervals are intervals with rational slope containing rational points, and intervals with irrational slope with one endpoint rational. If a conjecture by Franks and Misiurewicz is true \cite{fm}, then these are the only possible examples of rotation sets which are non-degenerate intervals.

Consider the simplest of such examples: a vertical interval of the form $\{0\}\times I$, containing $(0,0)$ in its interior. The canonical example of a homeomorphism with that rotation set is a twist $(x,y) \mapsto (x, y + \sin(2\pi x))$, which clearly has annular dynamics. One would like to see to what extent this is a model for a torus homeomorphism with a rotation set of the form $\{0\}\times I$. A partial answer to this question is given by Theorem \ref{teo.dav} below, which gives a qualitative description of the dynamics, in the case that $(0,0)$ is not realized by a periodic orbit. It proves in particular that the dynamics is annular. Before stating the theorem, we recall some definitions. We say that a curve $\gamma\subset\T^2$ is \textit{free forever} for $f$ if $f^n(\gamma)\cap\gamma=\empt$ for all $n\in\Z$. Also, if $\gamma_1,\gamma_2\subset\T^2$ are vertical and disjoint curves, using the covering $\R\times \T^1 \ra \T^2$ one can define the colsed annulus $[\gamma_1,\gamma_2]\subset\T^2$ whose `left' border component is $\gamma_1$ and whose `right' border component is $\gamma_2$. By last, denote by $\Omega(\tl{f})$ the \textit{non-wandering set} of $f$, that is, the set of points $x\in\T^2$ such that for every neighborhood $V$ of $x$, there is $n>0$ such that $f^n(V)\cap V\neq\empt$.

\begin{figure}[h]        
\begin{center} 
\includegraphics{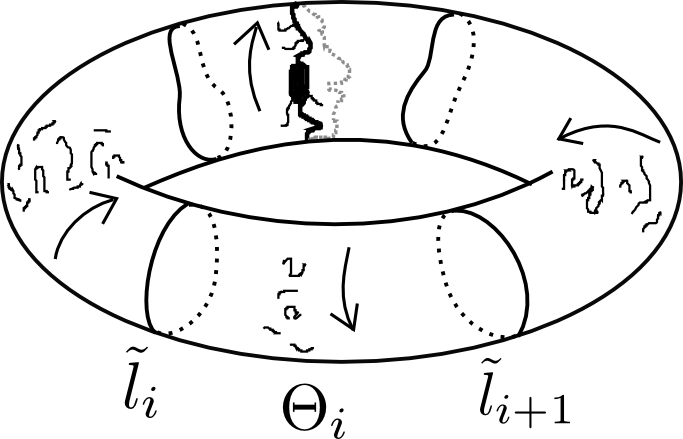}
\caption{Illustration for Theorem \ref{teo.dav}. At least one of the sets $\Theta_i$ must be annular and essential.}
\label{fig.teo1}
\end{center}  
\end{figure}

\begin{teorema}[\cite{dav}]    \label{teo.dav}
Let $f$ be a homeomorphism of $\T^2$ homotopic to the identity with a lift $\tl{f}:\R^2\ra\R^2$ such that:
\begin{itemize}
\item $\rho(\tl{f})=\{0\}\times I$, where $I$ is a non-degenerate interval containing $0$ in its interior, and
\item $(0,0)$ is not realized by a periodic point.
\end{itemize}
Then, the dynamics of $f$ is annular. Moreover, there exists a finite family $\{ l_i \}_{i=0}^{r-1}$, $r\geq 2$, of curves in $\T^2$ which are simple, closed, vertical, and pairwise dijoint, and with the following properties. If 
$$\Theta_i:= \bigcap_{n\in\Z} f^n \left( [l_i,l_{i+1}] \right)\ \ \ \txt{for } i\in\Z/r\Z,$$
then,
\begin{enumerate}
\item at least one of the sets $\Theta_i$ is an annular, essential, $f$-invariant set which is a semi-attractor, 
\item the curves $l_0,l_1,\ldots,l_{r-1}$ are free forever for $f$,
\item there is $\epsilon>0$ such that $\rho(\Theta_i,\tl{f})$ is contained either in $\{0\}\times (\epsilon,\infty)$, or in $\{0\}\times (-\infty,-\epsilon)$, and
\item $\Omega(f)\subset \cup \Theta_i$, (see Fig. \ref{fig.teo1}).
\end{enumerate}
\end{teorema}

This theorem tells us that, if $(0,0)\in\rho(\tl{f})$ is not realized by a periodic orbit, then there is a `filtration' for the dynamics of $f$, given by the curves $l_i$. At each level of the filtration, the maximal invariant set $\Theta_i$ rotates either upwards or downwards, and one of these maximal invariant sets is a `vertical barrier' for the dynamics. Also, in \cite{dav} is deduced from this an analogous theorem for the case that the rotation set is a general interval with rational slope and containing rationals.

In Theorem \ref{teo.dav}, the hypothesis that $(0,0)$ is not realized by a periodic orbit is equivalent to the fact that the lift $\tl{f}$ has no fixed points (this is Frank's Lemma \cite{f1}), which in turn is equivalent to the fact that $f$ has no contractible fixed points. Therefore, applying Theorem \ref{teo.pat} one easily obtains a foliation $\cl{F}$ of $\T^2$ such that the lift $\ell\subset \R^2$ of any leaf is a Brouwer curve for $\tl{f}$ (that is, $\tl{f}(\ell)\cap\ell= \tl{f}^{-1}(\ell)\cap\ell=\empt$). Then, in \cite{dav} it is proved that there exists a finite family of leaves $l_i$ of $\cl{F}$ that are compact, essential and vertical, and such that the maximal invariant sets $\Theta_i$ for $f$ between them rotate either upwards or downwards. This is done using the fact that $\cl{F}$ is positively transverse to the isotopy $I$ from Id to $f$, with techniques similar to those in \cite{lc3} and with the use of Atkinson's Lemma \cite{atk} from ergodic theory.

The main and new part in theorem \ref{teo.dav} is the existence of an annular essential $f$-invariant set. This is done in the following way. Let $\ell\subset\R^2$ be a lift of some of the curves $l_i\subset\T^2$. Then, $\ell$ is a free curve for $\tl{f}$. By contradiction, suppose that none of the sets $\Theta_i$ is an essential set. Then, one may easily see that 
$$\tl{f}^{n_0}(\ell)\cap T_1(\ell)\neq\empt,$$
for some $n_0\in\Z$ (assume without loss of generality that $n_0 >0$). The main work in \cite{dav} consists of proving that the fact that $\tl{f}^{n_0}(\ell)\cap T_1(\ell)\neq\empt$ implies actually that there exists $x\in\R^2$ such that 
$$ \lim_{n\ra\infty} \frac{\txt{pr}_1 (\tl{f}^n(x) - x)}{n} = \infty,$$
which yiels the contradiction $\max \txt{pr}_1(\rho(\tl{f})) > 0$. Hence, one of the sets $\Theta_i$ must be an annular essential set. The proof of the other items of Theorem \ref{teo.dav} follows from this.

For the general case of a homeomorphism with an annular rotation set, the following question is still open.

\begin{pregunta}    \label{preg1}
If a torus homeomorphism $f$ has an annular rotation set, then is the dynamics of $f$ annular?
\end{pregunta}

Progress in this direction has been made by Bortolatto an Tal \cite{tb}, showing that the answer to this question is affirmative in the case that Lebesgue measure is ergodic and has a non-rational rotation vector (that is, Lebesgue-almost every point in $\T^2$ has the same well defined non-rational rotation vector). More recently, Guelman, Koropecki and Tal have shown that the answer to Question \ref{preg1} is also affirmative if one assumes only that Lebesgue measure is preserved \cite{gkt}.

Now, for the case that the rotation set is an interval with irrational slope, as we mentioned above, the only known examples have a rational endpoint. Such an example can be given in the following way. Let $v\in\R^2$ be a vector with irrational slope, and let $\chi$ denote the constant vector field $\chi\equiv v$ in $\T^2$. For $p\in\T^2$, let $\psi:\T^2\ra\R$ be a continuous function such that $\psi \geq 0$ and $\psi(x)=0$ if and only if $x=p$. Now, let $f:\T^2\ra\T^2$ be the time-$1$ map of the flow given by the vector field $\psi\chi$. We have that $\txt{Fix}(f)=\{p\}$, and that the future orbit of every point passes arbitrarily close of $p$. Let $\tl{f}:\R^2\ra\R^2$ be the lift of $f$ wich fixes the lifts of $p$. If $\psi$ is chosen adequately (not too close to zero), one can prove that there are points in $\R^2$ with non-zero rotation vector for $\tl{f}$, and then the rotation set is an interval with irrational slope (equal to the slope of $v$) and with $(0,0)$ as an endpoint. 

In such example, one can of course replace the point $p$ by any totally disconnected closed set $X$, or even `explode' the set $X$ and some flux lines, remaining with examples with the same rotation set and which are semi-conjugate to the example constructed above. We don't know if there are any more examples with such a rotation set.

%\section{Rotation sets with interior}

%The only known examples of rotation sets with interior are usual polygons, and `polygons' with a countably infinite number of extremal points and sides. In \cite{kw1} it is proved that any polygon with rational vertices can be realized as a rotation set, and in \cite{kw2} it is constructed a torus homeomorphism whose rotation set is a `polygon' with a countably infinite number of rational extremal points.

%The canonical example of a torus homeomorphism with a polygon as a rotation set is a fitted Axiom A diffeomorphism. These type of diffeomorphisms admit always Markov partitions, which encode the rotation information. In \cite{rata} it is shown that any Axiom A diffeomorphism has always a rational polygon as rotation set. We observe that the example of an `infinite' poligon as a rotation set in \cite{kw2} is constructed using an infinite Markov partition, which also encodes the rotation information. 

%For the case that the rotation set has non-empty interior, we don't know wether there exists or not some dynamical information that may be considered as a dynamical model.

\bibliographystyle{amsalpha}
\bibliography{bib-cimat}

\end{document}